\documentclass[11pt,reqno]{amsart}

\usepackage[utf8]{inputenc}

\usepackage{amsmath,amsthm,amssymb}
\usepackage{amssymb}

\usepackage{graphics}
\usepackage{hyperref}
\usepackage[usenames, dvipsnames]{xcolor}

\usepackage{mathtools}
\mathtoolsset{showonlyrefs}

\usepackage[square,sort,comma,numbers]{natbib}

\definecolor{darkblue}{rgb}{0.0,0.0,0.3}
\hypersetup{colorlinks,breaklinks,
  linkcolor=darkblue,urlcolor=darkblue,
anchorcolor=darkblue,citecolor=darkblue}

\theoremstyle{plain}
\newtheorem{theorem}{Theorem}
\newtheorem*{theorem*}{Theorem}

\newtheorem*{proposition*}{Proposition}

\newtheorem*{corollary*}{Corollary}

\newtheorem*{openq*}{Open Question}
\newtheorem{conjecture}[theorem]{Conjecture}
\newtheorem*{conjecture*}{Conjecture}

\theoremstyle{definition}

\title{A Short Note on Gaps between Powers of Consecutive Primes}
\author{David Lowry-Duda}

\begin{document}

\maketitle

\begin{abstract}

  Let $\alpha, \beta \geq 0$ and $\alpha + \beta < 1$.
  In this short note, we show that
  $\liminf_{n \to \infty} p_n^\beta(p_{n+1}^\alpha - p_n^\alpha) = 0$,
  where $p_n$ is the $n$th prime.
  This notes an improvement over results of S\'{a}ndor and gives additional
  evidence towards a conjecture of Andrica.
  This follows directly from recent results on prime pairs from Maynard, Tao,
  Zhang.

\end{abstract}

\section{Introduction}

The primary purpose of this note is to collect a few hitherto unnoticed or
unpublished results concerning gaps between powers of consecutive primes.
The study of gaps between primes has attracted many mathematicians and led to
many deep realizations in number theory.
The literature is full of conjectures, both open and closed, concerning the
nature of primes.

In a series of stunning developments, Zhang, Maynard, and
Tao~\cite{Maynard, Zhang} made the first major progress towards proving the
prime $k$-tuple conjecture, and successfully proved the existence of infinitely
many pairs of primes differing by a fixed number.
As of now, the best known result is due to the massive collaborative Polymath8
project~\cite{Polymath8}, which showed that there are infinitely many pairs of
primes of the form $p, p+246$.
In the excellent expository article~\cite{Granville}, Granville describes the
history and ideas leading to this breakthrough, and also discusses some of the
potential impact of the results.
This note should be thought of as a few more results following from the
ideas of Zhang, Maynard, Tao, and the Polymath8 project.

Throughout, $p_n$ will refer to the $n$th prime number.
In~\cite{Andrica}, Andrica conjectured that
\begin{equation}\label{eq:Andrica_conj}
  \sqrt{p_{n+1}} - \sqrt{p_n} < 1
\end{equation}
holds for all $n$.
This conjecture, and related statements, is described in~\cite{GuyUnsolved}.
It is quickly checked that this holds for primes up to $4.26 \cdot 10^{8}$
using sagemath~\cite{sagemath}.\footnote{Code verifying this is available on the
author's website at \url{http://davidlowryduda.com/?p=2430}}
It appears very likely that the conjecture is true.
However it is also likely that new, novel ideas are necessary before the
conjecture is decided.

Andrica's Conjecture can also be stated in terms of prime gaps.
Let $g_n = p_{n+1} - p_n$ be the gap between the $n$th prime and the $(n+1)$st
prime.
Then Andrica's Conjecture is equivalent to the claim that
$g_n < 2 \sqrt{p_n} + 1$.
In this direction, the best known result is due to Baker, Harman, and
Pintz~\cite{BakerHarmanPintz}, who show that $g_n \ll p_n^{0.525}$.

In 1985, S\'{a}ndor~\cite{sandor1985certain} proved that
\begin{equation}\label{eq:Sandor}
  \liminf_{n \to \infty} \sqrt[4]{p_n} (\sqrt{p_{n+1}} - \sqrt{p_n}) = 0.
\end{equation}
The close relation to Andrica's Conjecture~\eqref{eq:Andrica_conj} is clear.
The first result of this note is to strengthen this result.

\begin{theorem}
  Let $\alpha, \beta \geq 0$, and $\alpha + \beta < 1$.
  Then
  \begin{equation}\label{eq:main}
    \liminf_{n \to \infty} p_n^\beta (p_{n+1}^\alpha - p_n^\alpha) = 0.
  \end{equation}
\end{theorem}

We prove this theorem in Section~\ref{sec:proof}.
Choosing $\alpha = \frac{1}{2}, \beta = \frac{1}{4}$ verifies S\'{a}ndor's
result~\eqref{eq:Sandor}.
But choosing $\alpha = \frac{1}{2}, \beta = \frac{1}{2} - \epsilon$ for a small
$\epsilon > 0$ gives stronger results.

This theorem leads naturally to the following conjecture.

\begin{conjecture}\label{conj}
  For any $0 \leq \alpha < 1$, there exists a constant $C(\alpha)$ such that
  \begin{equation}
    p_{n+1}^\alpha - p_{n}^\alpha \leq C(\alpha)
  \end{equation}
  for all $n$.
\end{conjecture}

A simple heuristic argument, given in Section~\ref{sec:heuristic}, shows that
this Conjecture follows from Cram\'{e}r's Conjecture.

It is interesting to note that there are generalizations of Andrica's
Conjecture.
One can ask what the smallest $\gamma$ is such that
\begin{equation}
  p_{n+1}^{\gamma} - p_n^{\gamma} = 1
\end{equation}
has a solution.
This is known as the Smarandache Conjecture, and it is believed that the
smallest such $\gamma$ is approximately
\begin{equation}
  \gamma \approx 0.5671481302539\ldots
\end{equation}
The digits of this constant, sometimes called ``the Smarandache constant,'' are
the contents of sequence A038458 on the OEIS~\cite{OEIS}.
It is possible to generalize this question as well.

\begin{openq*}
  For any fixed constant $C$, what is the smallest $\alpha = \alpha(C)$ such that
  \begin{equation}
    p_{n+1}^\alpha - p_n^\alpha = C
  \end{equation}
  has solutions?
  In particular, how does $\alpha(C)$ behave as a function of $C$?
\end{openq*}

This question does not seem to have been approached in any sort of generality,
aside from the case when $C = 1$.

\section{Proof of Theorem}\label{sec:proof}

The idea of the proof is very straightforward.
We estimate~\eqref{eq:main} across prime pairs $p, p+246$, relying on the recent
proof~\cite{Polymath8} that infinitely many such primes exist.

Fix $\alpha, \beta \geq 0$ with $\alpha + \beta < 1$.
Applying the mean value theorem of calculus on the function $x \mapsto x^\alpha$
shows that
\begin{align}
  p^\beta \big( (p+246)^\alpha - p^\alpha \big)
  &=
  p^\beta \cdot 246 \alpha q^{\alpha - 1} \nonumber
  \\
  &\leq
  p^\beta \cdot 246 \alpha p^{\alpha - 1}
  =
  246 \alpha p^{\alpha + \beta - 1}, \label{eq:bound}
\end{align}
for some $q \in [p, p+246]$.
Passing to the inequality in the second line is done by realizing that
$q^{\alpha - 1}$ is a decreasing function in $q$.
As $\alpha + \beta - 1 < 0$, as $p \to \infty$ we see that~\eqref{eq:bound} goes
to zero.

Therefore
\begin{equation}
  \liminf_{n \to \infty} p_n^\beta (p_{n+1}^\alpha - p_n^\alpha) = 0,
\end{equation}
as was to be proved.

\section{Further Heuristics}\label{sec:heuristic}

Cram\'{e}r's Conjecture states that there exists a constant $C$ such that for
all sufficiently large $n$,
\begin{equation}
  p_{n+1} - p_n < C(\log n)^2.
\end{equation}
Thus for a sufficiently large prime $p$, the subsequent prime is at most
$p + C (\log p)^2$.
Performing a similar estimation as in Section~\ref{sec:proof} shows that
\begin{equation}
  (p + C (\log p)^2)^\alpha - p^\alpha \leq C (\log p)^2 \alpha p^{\alpha - 1} =
  C \alpha \frac{(\log p)^2}{p^{1 - \alpha}}.
\end{equation}
As the right hand side vanishes as $p \to \infty$, we see that it is natural to
expect that Conjecture~\ref{conj} is true.
More generally, we should expect the following, stronger conjecture.
\begin{conjecture}
  For any $\alpha, \beta \geq 0$ with $\alpha + \beta < 1$, there exists a
  constant $C(\alpha, \beta)$ such that
  \begin{equation}
    p_n^\beta (p_{n+1}^\alpha - p_n^\alpha) \leq C(\alpha, \beta).
  \end{equation}
\end{conjecture}

\vspace{20 mm}
\bibliographystyle{alpha}
\bibliography{shortbib}

\end{document}